\documentclass[twoside,a4paper]{amsart}

\PassOptionsToPackage{pdfauthor={Vladimir V. Kisil},%
    pdftitle={Metamorphism---an Integral Transform Reducing the Order of a Differential Equation},%
    pdfsubject={mathematics},%
    backref=page,%
  pdfkeywords={symmetries}
}{hyperref}

\usepackage[utf8]{inputenc}				
\newif\ifdgrutter
\ifdgrutter
\usepackage[big,online]{dgruyter}	
\fi
\usepackage{lmodern} 
\usepackage{microtype}
\usepackage[numbers,square,sort&compress]{natbib}

\IfFileExists{eulervm.sty}{%
\usepackage{eulervm}}{}
\usepackage{hyperref,graphicx,colonequals}

\numberwithin{equation}{section}

\providecommand{\gpartial}[2]{\partial_{\ifcase#2
  x \or y \or b \or r \fi}}

\providecommand{\half}{{\textstyle\frac{1}{2}}}
\providecommand{\third}{{\textstyle\frac{1}{3}}}
   \PassOptionsToPackage{british}{babel}
\theoremstyle{dgthm}
\newtheorem{thm}{Theorem}
    \newtheorem{prop}[thm]{Proposition}
    \newtheorem{lem}[thm]{Lemma}
    \newtheorem{cor}[thm]{Corollary}
\theoremstyle{dgdef}

    \newtheorem{example}[thm]{Example}

   \theoremstyle{remark}
    \newtheorem{rem}[thm]{Remark}
    
\providecommand{\norm}[2][\relax]{\left\|#2\right\|\ifx#1\relax\else_{#1}\fi}
\providecommand{\modulus}[2][\relax]{\left| #2 \right|\ifx#1\relax\else_{#1}\fi}
\providecommand{\oper}[1]{\mathcal{#1}}

\providecommand{\Space}[3][]{\ifx#2R\ifx#1e \mathbb{C}^{#3} \else
\ifx#1p \mathbb{D}^{#3} \else
\ifx#1h \mathbb{O}^{#3} \else
\ifx#1\sigma \mathbb{A}\!^{#3} \else
\ensuremath{\mathbb{#2}^{#3}_{#1}{}} \fi \fi \fi \fi \else
\ensuremath{\mathbb{#2}^{#3}_{#1}{}} \fi}
\providecommand{\FSpace}[3][]{\ensuremath{\ifx#2l \ell_{#3}^{#1}{}\else
  \mathsf{#2}_{#3}^{#1}{}\fi}}

\providecommand{\uir}[3][0]{\ifcase #1{\rho^{#2}_{#3}}%
\or {\breve{\rho}^{#2}_{#3}}%
\or {\tilde{\rho}^{#2}_{#3}}\fi}
\usepackage{braket}
\providecommand{\scalar}[3][\relax]{\left\langle #2, #3
        \right\rangle\ifx#1\relax\else_{#1}\fi}

\providecommand{\rme}{\mathrm{e}}
\providecommand{\rmd}{\mathrm{d}}

\providecommand{\myhbar}{\hbar}
\newcommand{\kbar}{
  {\mathchar'26}\mkern-7mu k}
\makeatletter
\newcommand{\@hslashslash}{%
  \raisebox{-0.2ex}{%
    \scalebox{1}{%
      \rotatebox[origin=c]{20}{$\mathchar'26$%
      }%
    }%
  }%
}
\newcommand{\@kslash}[2]{%
  {%
   \vphantom{#2}%
   \ooalign{\kern#1em\smash{\@hslashslash}\hidewidth\cr
     $#2$\cr
   }%
   \kern.05em
  }%
}
\newcommand\kslash{\mathchoice
  {\hbox{\@kslash{.01}{k}}}
  {\hbox{\@kslash{.01}{k}}}
  {\hbox{\fontsize{\sf@size}{\sf@size}\selectfont\@kslash{.05}{k}}}
  {\hbox{\fontsize{\ssf@size}{\ssf@size}\selectfont\@kslash{.1}{k}}}
}
\let\kbar=\kslash

\providecommand{\myeprint}[2]{\href{#1}{\texttt{#2}}}
\DeclareFontFamily{OT1}{cyr}{}
\DeclareFontShape{OT1}{cyr}{m}{n}
   {  <5> <6> <7> <8> <9> gen * wncyr
      <10> <10.95> <12> <14.4> <17.28> <20.74> <24.88> wncyr10}{}
\DeclareFontShape{OT1}{cyr}{m}{it}
    {
       <5> <6> <7> <8> <9> gen * wncyi
      <10> <10.95> <12> <14.4> <17.28> <20.74> <24.88>wncyi10
      }{}
\DeclareFontShape{OT1}{cyr}{m}{ss}
    {
       <5> <6> <7> <8> wncyss8
       <9> wncy9
      <10> <10.95> <12> <14.4> <17.28> <20.74> <24.88>wncyss10
      }{}
\DeclareFontShape{OT1}{cyr}{m}{sc}
    {
       <5> <6> <7> <8> <9> <10> <10.95> <12> <14.4> <17.28> <20.74> <24.88>wncysc10
      }{}
\DeclareFontShape{OT1}{cyr}{bx}{n}
   {
       <5> <6> <7> <8> <9> gen * wncyb
      <10> <10.95> <12> <14.4> <17.28> <20.74> <24.88>wncyb10
      }{}
\DeclareTextFontCommand{\textcyr}{\fontfamily{cyr}\selectfont}

\providecommand{\cprime}{'}

\usepackage{accents}
\newcommand{\metamorph}[1]{\accentset{\smash{\raisebox{-0.12ex}{$\scriptscriptstyle\approx$}}}{#1}\rule{0pt}{2.3ex}}

\providecommand{\doublechecked}{}

\let\pyginac=\iffalse
\let\endpyginac=\fi
\providecommand{\cites}[1]{\cite{#1}}
\providecommand{\citelist}{}
\providecommand{\amscite}[3]{\cite[#3]{#1}}
\begin{document}
\ifdgrutter

	\articletype{Research Article}
	\received{Month	DD, YYYY}
	\revised{Month	DD, YYYY}
  \accepted{Month	DD, YYYY}
  \journalname{De~Gruyter~Journal}
  \journalyear{YYYY}
  \journalvolume{XX}
  \journalissue{X}
  \startpage{1}
  \aop
  \DOI{10.1515/sample-YYYY-XXXX}

  \title{Metamorphism---an Integral Transform\\ Reducing the Order of a Differential Equation}
\runningtitle{Metamorphism and Helmholtz Equation}

\author[1]{Vladimir V. Kisil}
\runningauthor{V.V. Kisil}

\affil[1]{\protect\raggedright 
School of Mathematics,
University of Leeds,
Leeds LS2\,9JT,
England, email: {\href{mailto:V.Kisil@leeds.ac.uk}{V.Kisil@leeds.ac.uk}}, 
{\href{mailto:kisilv@maths.leeds.ac.uk}{kisilv@maths.leeds.ac.uk}}
}
\abstract{  We propose an integral transform, called metamorphism, which allows us to reduce the order of a differential equation. For example, the second order Helmholtz equation is transformed into a first order equation, which can be solved by the method of characteristics.}
\else
  \title[Metamorphism and Helmholtz Equation]{Metamorphism---an Integral Transform\\ Reducing the Order of a Differential Equation}

\author[Vladimir V. Kisil]{\href{http://www1.maths.leeds.ac.uk/~kisilv/}{Vladimir V. Kisil}}

\address{
School of Mathematics,
University of Leeds,
Leeds LS2\,9JT,
England}

\email{\href{mailto:V.Kisil@leeds.ac.uk}{V.Kisil@leeds.ac.uk}}
\urladdr{\url{http://www1.maths.leeds.ac.uk/~kisilv/}}

\date{\today}

\begin{abstract}
We propose an integral transform, called metamorphism, which allows us to reduce the order of a differential equation. For example, the second order Helmholtz equation is transformed into a first order equation, which can be solved by the method of characteristics.
\end{abstract}
\subjclass{Primary: 35A22; Secondary: 35C15}
\fi
\keywords{Partial differential equation, integral transform, Helmholtz equation, coherent states}

\maketitle


\section{Introduction}
\label{sec:introduction}

Partial differential equations (PDEs) provide a fundamental language describing laws of nature. To tackle intrinsic complexity of PDEs one often tries to transform an equation to a simpler one. For example:
\begin{enumerate}
\item \emph{Separation of variables}~\cite{Miller77a} allows one to replace a given PDE by several other differential equations, each with a smaller number of variables. Ideally, one wants to obtain a system of ordinary differential equations (ODEs), which are much more accessible.
\item \emph{Integral transformations}~\amscite{PolyaninNazaikinskii16a}*{\S~15.2} which map some derivatives to simpler objects, e.g. multiplication by variables. The Fourier and Laplace transforms are the most common choice, with the Mellin, Hankel, etc. transforms to follow in more specialised cases.
\item \emph{Special functions}~\cites{Miller68,Vilenkin68} (including \emph{orthogonal polynomials}) are a blend of the two previous techniques. A decomposition over a suitable family of special functions is a sort of integral transform and level curves of the special functions define separating coordinates for PDEs, cf.~\cite{ColbrookKisil20a}.
\item \emph{Transmutations}  relate solutions of a differential equation with variable coefficients  to solutions of an equation with constant coefficients. The latter admits a better understanding, thus such a reduction is very helpful starting from second order ODEs~\cite{KravchenkoSitnik20a}. 
\end{enumerate}
This paper adds to this list a generic method, called \emph{metamorphism}, to study PDEs. It is an integral transform~\eqref{eq:metamorphism} which encompasses the power of many familiar maps, see Rem.~\ref{re:generalisations} below.
 Metamorphism allows one to reduce the order of PDEs. For example, a second order PDE can be reduced to a first order one, which opens doors for other techniques, e.g. the method of characteristics. The approach is illustrated by application to the Helmholtz equation.

\section{Metamorphism}
\label{sec:metamorphism}

The key ingredient of our method is an integral transform called \emph{metamorphism}. It is a covariant transform (aka coherent state transform)~\cite{AliAntGaz14a} for a certain subgroup of the Jacobi--Schr\"odinger group~\citelist{\amscite{Folland89}*{\S~4.4} \amscite{Berndt07a}*{\S~8.5}}. This was implemented in the working Jupyter notebook~\cite{Kisil21b} and is spelled in~\cite{AlqurashiKisil23a}. Here, all required properties of the metamorphism will be verified by direct arguments without a reference to the representation theory.

\subsection{The Integral Transform}
\label{sec:integral-transform}

In the following presentation the parameter \(\myhbar\) shall be associated with the wave length (or the Planck constant) and it is a fixed parameter most of the time. 
The \emph{metamorphism} is an integral operator defined by:
\begin{equation}
  \label{eq:metamorphism}
  \begin{split}
    \metamorph
    {f}(x,y,b,r)
  &={\sqrt{\myhbar r}}
  \int_{\Space{R}{}} f(u)\,
  \exp\left( -\pi  \myhbar \left(
  ( r^{2}-i b)(u-y)^2
   +2 i  (u-y) x
  \right)\right)
  \,\rmd u
  \,.
\end{split}
\end{equation}
Here \(x\), \(y\), \(b\) are reals and \(r\) is a positive real, their collection is denoted by \(\Space[+]{R}{4}\).  The integral~\eqref{eq:metamorphism} is meaningful for \(f(u)\) from many linear spaces, e.g. the Schwartz space,   \(\FSpace{L}{p}(\Space{R}{})\) with \(1\leq p\leq \infty\). Being based on the Gaussian, metamorphism~\eqref{eq:metamorphism} is also well-defined for tempered distributions from \(\FSpace[\prime]{S}{}(\Space{R}{})\) and maps them to smooth functions.  This transformation is a covariant (aka coherent states) transform~\cite{AliAntGaz14a}, namely it has the form
\begin{equation}
  \label{eq:metamorphism-scalar}
    \metamorph
    {f}(x,y,b,r)= \scalar{f}{ \phi_{x,y,b,r}}, \qquad
 \text{ for }\quad \phi_{x,y,b,r}= \uir{}{x,y,b,r} \phi\,,
\end{equation}
where \(\phi(u)=\rme^{-\pi\myhbar u^2}\) is the Gaussian and \(\uir{}{x,y,b,r}\) is a unitary irreducible representation of the SSR group~\cite{Kisil21b}---the semidirect product of the Heisenberg group \(\Space[p]{H}{1}\) with the \(ax+b\) group acting by symplectic automorphisms of \(\Space[p]{H}{1}\)~\amscite{Folland89}*{\S~1.2}.
Alternatively, the metamorphism can be  stated as the convolution
\begin{displaymath}
  \metamorph
    {f}(x,y,b,r)
  = [f * K_{x,b,r}](y)
\end{displaymath}
with the Gauss-type kernel
\begin{displaymath}
  K_{x,b,r}(u) = {\sqrt{r}}\,
  \rme^{-\pi  \myhbar (
  (r^{2}-i b)u^2
   -2 i  u x
   )}
 =[\uir{}{x,0,-b,r} \phi](u)\,
\end{displaymath}
parametrised by \(x\), \(b\in\Space{R}{}\) and \(r\in\Space[+]{R}{}\). Note the reverse sign of \(b\) in the last formula.


\begin{rem}
  \label{re:generalisations}
Metamorphism~\eqref{eq:metamorphism} incorporates many known integral transformations:
\begin{itemize}
\item For \(r\rightarrow 0\) and \(y=b=0\) it reduces to the Fourier transform.
\item For \(x=b=0\) we get the Gauss--Weierstrass(--Hille) transform~\cite{Zemanian67a} which solves the Cauchy problem for the heat equation. Its imaginary cousins with \(x=0\) and \(r \rightarrow 0\) solve the time evolution of a free quantum particle~\amscite{Miller77a}*{Ch.1~(1.25)}.
\item For \(r=1\), \(b=0\) this is the Fock--Segal--Bargmann (FSB) transform~\cites{Segal60,Bargmann61}, see also~\citelist{\amscite{Folland89}*{\S~1.6}}.
\item For \(b=0\) (with variable \(r\)) it is the Fourier--Bros--Iagolnitzer (FBI) transform~\amscite{Folland89}*{\S~3.3}.
\item For \(b=0\) and \(r = (1+\modulus{x})^\alpha\) for some \(\alpha\) the metamorphism asymptotics for  \(x \rightarrow \infty\) is equivalent to the C\'{o}rdoba--Fefferman treatment of wavefronts~\citelist{\cite{CordobaFefferman78a} \amscite{Folland89}*{\S~3.1}}.  
\item With \(r=1\) it was considered in~\cites{AlmalkiKisil18a,AlmalkiKisil19a} to treat the Schr\"odinger equation.
\item The metamorphism is also connected to the wavelet transform for the affine group~\amscite{AliAntGaz14a}*{Ch.~12}, which is a subgroup of the SSR group. 
\end{itemize}
Integral transforms of this type were extensively studied~\cites{Neretin11a,Osipov92a} and applied in many areas~\cites{HealyKutayOzaktasSheridan16,Folland89}. Although metamorphism can be considered as a special type of a complex linear canonical transform (LCT)~\cite{HealyKutayOzaktasSheridan16}, 
the particular form~\eqref{eq:metamorphism} still has a large unexplored potential in applications.
\end{rem}

Metamorphisms can be computed explicitly for some important functions.

\begin{example}
  \label{ex:initial-metamorphisms}
  \begin{enumerate}
  \item We start from the wave packet \(P_{\sigma,\lambda}(u)= \rme^{-\pi \sigma u^2 -2\pi i \lambda u} \in \FSpace{L}{2}(\Space{R}{})\) for a complex \(\sigma\) such that \(\Re \sigma>0\):
    \doublechecked
    \begin{equation}
\pyginac
sigma=realsymbol("sigma")
u=realsymbol("u")
h=realsymbol("h")
y=realsymbol("y")
la=realsymbol("lambda")
W=-Pi*sigma*u**2-2*Pi*I*la*u
E=-Pi*sigma*(u-y)**2-2*Pi*h*I*(u-y)*(la-I*sigma*y)/h
print((W-E).expand().normal())
# - \pi y^{2} \sigma-{(2 i)}  \pi y \lambda
\endpyginac
      \label{eq:metamorphism-wave-packet}
      \begin{split}
      \metamorph{P}_{\sigma,\lambda}  (x,y,b,r)
      & = \frac{ \sqrt{r}}{\sqrt{r^{2}+\sigma/\myhbar-i b}}
        \exp(-  \pi \sigma  y^2 - 2 \pi i \lambda y)\\
      & \qquad \times 
      \exp\left(- \frac{\pi\myhbar (x +(\lambda-i\sigma y)/\myhbar)^{2}}{r^{2}+\sigma/\myhbar-i b}\right)\\
      &=
      \frac{\sqrt{ r}}{\sqrt{ r^{2}+\sigma / \myhbar -i  b}}
      \exp(  \pi i  \myhbar  (2   x y + (b+i r^{2}) y^2))\\
      &\qquad \times
       \exp\left(- \frac{  \pi \myhbar {(x+ (b + i r^{2} ) y- i \lambda/\myhbar)}^{2}}{ r^{2}+\sigma/\myhbar  -i    b }\right)\,.
    \end{split}
  \end{equation}
  This follows from the well-known formula~\amscite{Folland89}*{App.~A, (1)}:
  \begin{displaymath}
    \int_{\Space{R}{}} \rme^{-\pi a u^2 - 2 \pi i z u}\,\rmd u = \frac{1}{\sqrt{a}}\, \rme^{-\pi z^2/a}\,,
  \end{displaymath}
  where \(\Re a> 0\) and \(z\in\Space{C}{}\) with the usual agreement on the branch of \(\sqrt{a}\).
  \item For the exponent \(E_{k}(u)=\rme^{- i k u}=\rme^{-2 \pi i \kbar u}\) (with \(\kbar=k/(2 \pi )\)) representing a wave with wave number \(k\) we have:
    \doublechecked
    \begin{align}
      \label{eq:wave-transformed}
      \metamorph{E}_k (x,y,b,r) &= \frac{ \sqrt{r}}{\sqrt{r^{2}-i b}}  \exp(- i k  y) \exp(- \pi\myhbar (k/(2 \pi \myhbar) +x)^{2}/(r^{2}-i b))\\
      \nonumber 
      &= \frac{ \sqrt{r}}{\sqrt{r^{2} - i b}}  \exp(- 2 \pi i \kbar  y) \exp(- \pi\myhbar (\kbar/ \myhbar +x)^{2}/(r^{2}-i b))\,. 
    \end{align}
    It formally coincides with the transform of the wave packet \(P_{\sigma,\lambda}\) with \(\sigma=0\) and \(\lambda = \kbar\).

  \item
    For the Dirac delta function \(\delta\) and its derivative \(\delta_1\) we have:
    \begin{align*}
      \metamorph{\delta} (x,y,b,r)
      &= \sqrt{r}
        { \rme^{-\pi \myhbar (r^2 - i b )  y^{2}} 
        \rme^{{2 \pi i \myhbar }  y   x}}\,,\\
      \metamorph{\delta}_1 (x,y,b,r)
      &=
        -2  \pi i  \myhbar \sqrt{r}{\left(x + (b+i r^2) y\right)}
        { \rme^{-\pi \myhbar (r^2 - i b )  y^{2}} 
    \rme^{{2 \pi i \myhbar }  y   x}}\,.
    \end{align*}
  \end{enumerate}
It is shown in \S~\ref{sec:char-image-space} that the repeated pattern \(x + (b+i r^2) y\) in the above formulae is not accidental.
\end{example}

\subsection{Sesqui-unitarity and inverse metamorphism}
\label{sec:sesqui-unitarity-fsb}
The known sesqui-unitarity property of Fourier--Wigner transform~\amscite{Folland89}*{\S~I.4} implies that:
\pyginac
def Gauss_integral(E,u):
    E1 = E[1].expand()
    a = -coeff(E1, u, 2)/Pi
    z = -coeff(E1, u, 1)/(2*Pi*I)
    c = coeff(E1, u, 0)
    
    return [E[0]/sqrt(a), c - Pi*z**2/a]

s=realsymbol("s","s")
x=realsymbol("x","x")
y=realsymbol("y", "y")
b=realsymbol("b", "b")
r=possymbol("r", "r")

half = numeric(1,2)
h=realsymbol("h","\\myhbar")
    
Metamorphism = [sqrt(r/(r**2+half-I*b)),\
(-Pi*h*(x+(b+I*r**2)*y)**2/(r**2+half-I*b)\
+Pi*I*h*(2*y*x+(b+I*r**2)*y**2))]

Metamorphism2 = [(r/abs(r**2+half-I*b)),\
2*(-Pi*h*(x+(b+I*r**2)*y)**2/(r**2+half-I*b)\
+Pi*I*h*(2*y*x+(b+I*r**2)*y**2))]

Int_y = Gauss_integral(Metamorphism, y)
Int_xy = Gauss_integral(Int_y, x)

print("

# => (-(r**2*h-I*b*h)**(-1)*(-1+(2*I)*b-2*r**2))**(1/2)*2**(1/2)*(-r*(-1+(2*I)*b-2*r**2)**(-1))**(1/2)*((-I*b+r**2)*h**(-1))**(1/2) * exp(0)
# = (2 r)^(1/2) / h
\endpyginac
\begin{align}
  \nonumber 
  \scalar{f_1}{f_2}
  & = \int_{\Space{R}{2}}
    \metamorph {f}_1(x,y,b_0,r_0)\,
    \overline{\metamorph {f}_2(x,y,b_0,r_0)}\, \frac{\myhbar\, \rmd x \,\rmd y}{\sqrt{2 r_0}}\,,\\
  \label{eq:sesqui-linearity-FW}
  & \equalscolon
    \scalar[(b_0,r_0)]
    {
    \metamorph {f}_1}
    {\metamorph {f}_2}\,,
    \quad \text{for any fixed } (b_0,r_0)\in\Space[+]{R}{2} \colonequals \Space{R}{}\times \Space[+]{R}{}\,. 
\end{align}
As an immediate consequence we obtain the (phase-space) \emph{reproducing property} for the metamorphosis, cf~\eqref{eq:metamorphism-scalar}:
\begin{align}
  \nonumber 
  \metamorph {f}(x,y,b,r)
  & = \scalar[\Space{R}{}]{f}{\phi_{x,y,b,r}}\\
  \label{eq:metamorphism-reproducing}
  & = \scalar[(b_0,r_0)]
    {\metamorph{f}}{\metamorph{\phi}_{x,y,b,r}}
\end{align}
with the reproducing kernel
\begin{align*}
  \metamorph{\phi}_{x,y,b,r}(x_0,y_0,b_0,r_0)
  & = \scalar[\Space{R}{}]{\phi_{x,y,b,r}}{\phi_{x_0,y_0,b_0,r_0}} = \overline{\metamorph{\phi}_{x_0,y_0,b_0,r_0}(x,y,b,r)}\,.
\end{align*}


\pyginac
def Gauss_integral(E,u):
    E1 = E.expand()
    a = -coeff(E1, u, 2)/Pi
    z = -coeff(E1, u, 1)/(2*Pi*I)
    c = coeff(E1, u, 0)
    
    return 1/sqrt(a)*exp(c)* exp (-Pi*z**2/a)

def repres(s, x, y, b, r, u):
    return -2*Pi*I*h*s-Pi*h*((r**2+I*b)*(u-y)**2-2*I*(u-y)*x)

sigma=realsymbol("sigma")
u=realsymbol("u")
h=realsymbol("h")

x=realsymbol("x")
y=realsymbol("y")
r=realsymbol("r")
b=realsymbol("b")

s1=realsymbol("s'","s")
x1=realsymbol("x'","x")
y1=realsymbol("y'", "y")
b1=realsymbol("b'", "b")
r1=possymbol("r'", "r")

s2=realsymbol("s_2")
x2=realsymbol("x_2")
y2=realsymbol("y_2")
b2=realsymbol("b_2")
r2=possymbol("r_2")

s0=realsymbol("s_0")
x0=realsymbol("x_0")
y0=realsymbol("y_0")
b0=realsymbol("b_0")
r0=possymbol("r_0")

z=symbol("z")
w=symbol("w")
z0=symbol("z0", "z_0")
w0=symbol("w0", "w_0")
la=realsymbol("lambda")
s=realsymbol(s)

# We substitute the parameters for the coherent state kernel
# and the composition of tbwo points
P=-Pi*h*u**2-2*Pi*I*h*s-Pi*h*((r**2-I*b)*(u-y)**2+2*I*(u-y)*x)

print(Gauss_integral(P,u).normal())
# We sunstitute the parameters of the left action on the SSR group
P_prod=Gauss_integral(P,u)\
.subs({s : y2**2*b1/2-y2*y1*b1-y2*x1+x1*y1+y1**2*b1/2,\
 x : -r1**(-1)*(x1+y1*b1-x2-y2*b1), y : r1*y2-r1*y1, b : r1**(-2)*(b2-b1), r : r1**(-1)*r2 })

print("Substituted point: $

## Manualy entered form of the reproducing kernel
Answer = exp(Pi*h*(x2-x1+(b2+I*r2**2)*(y2-y1))**2/(r2**2+r1**2+I*(b1-b2)))\
*exp(-Pi*I*h*(2*(y2-y1)*x2+(b2+I*r2**2)*(y2-y1)**2))

print("Comaprison 1: 
exp(wild(3))*wild(5)*exp(wild(4))==exp(wild(3)+wild(4))*wild(5)],subs_options.algebraic)\
.subs([exp(wild(1))*exp(wild(2))==exp(wild(1)+wild(2)),\
exp(wild(3))*wild(5)*exp(wild(4))==exp(wild(3)+wild(4))*wild(5)],subs_options.algebraic).expand().normal())

## Check the second source of the reproducing kernel
P_prod2 = Gauss_integral(repres(0, x1, y1, b1, r1, u)+conjugate(repres(0, x2, y2, b2, r2, u)), u)
math_filter(P_prod2.subs(h==1).normal())

print("Comaprison 2: 
exp(wild(3))*wild(5)*exp(wild(4))==exp(wild(3)+wild(4))*wild(5)],subs_options.algebraic)\
.subs([exp(wild(1))*exp(wild(2))==exp(wild(1)+wild(2)),\
exp(wild(3))*wild(5)*exp(wild(4))==exp(wild(3)+wild(4))*wild(5)],subs_options.algebraic).expand().normal())

\endpyginac

Utilising formula~\eqref{eq:metamorphism-wave-packet} with \(\sigma=\myhbar(r_0^2+i b_0)\) and \(\lambda = -\myhbar(x_0+(b_0-ir_0^2)y_0)\) we get:
\doublechecked
\begin{displaymath}
  \begin{split}
  \metamorph{\phi}_{x,y,b,r}(x_0,y_0,b_0,r_0)
& = \sqrt{\frac{r r_0}{ r_0^{2} +  r^{2}+i   (b-  b_0)}}\\
  & \qquad \times  \exp\left(-\frac{   \pi \myhbar {( x_0 - x +  (b_0 + i r_0^{2}) (y_0-  y))}^{2}}{
    r_0^{2}+r^{2}+i (b-b_0)}\right)\\
  & \qquad \times \exp\left(\pi i \myhbar \left(
    {2 }  (y_0-y)   x_0
    +(  b_0+i  r_0^{2}) ( y_0- y)^{2}
  \right)\right).
\end{split}
\end{displaymath}

The significance of~\eqref{eq:metamorphism-reproducing} is that we are able to restore \(\metamorph {f}(x,y,b,r)\) for any values of \((b,r)\) if \(\metamorph {f}(x,y,b_0,r_0)\) is initially known only for particular \((b_0,r_0)\) but any \((x,y)\).

The metamorphism \(\metamorph{f}\) in~\eqref{eq:metamorphism} contains an abundance of information on the function \(f\). Therefore, a function can be recovered from its metamorphism in many different ways. For example, we can work from the reconstruction formula for the FSB space: sesqui-unitarity~\eqref{eq:sesqui-linearity-FW} implies that the inverse metamorphism is provided by the adjoint operator for the pairing \(\scalar[(b_0,r_0)]{\cdot}{\cdot}\). Thus, the original function \(f(u)\) can be recovered from its transformation \(\metamorph{f}(x,y,b_0,r_0)\) for arbitrary fixed values \(b_0\in\Space{R}{}\) and \(r_0\in\Space[+]{R}{}\) as follows:
\begin{align}
  \label{eq:metamorphism-inverse}
  f(u)
  &=
  {\sqrt{r}}
  \int_{\Space{R}{2}} \metamorph{f}(x,y,b_0,r_0)\,
  \rme^{-\pi  \myhbar (
  (r_0^{2}+i b_0)(u-y)^2
   -2 i  (u-y) x
 )}
  \,\rmd x \,\rmd y
\end{align}
with integration over the phase space.

\subsection{Characterisation of the image space}
\label{sec:char-image-space}

The metamorphism integral kernel and thus the image space is annihilated by the following differential operators:
\begin{align}
  \label{eq:Cauchy-first}
  \oper{C}_1  & = \frac{1}{r}\left({(r^{2}-i b)} \gpartial_{0}+i \gpartial_{1}+2  x \myhbar \pi I\right) ;\\
  \label{eq:Cauchy-second}
  \oper{C}_2 &=  2  r^{2} \gpartial_{2} +i  r \gpartial_{3} - \half i I\,. 
\end{align}
 It is convenient to view these operators as the Cauchy--Riemann-type operators for complex variables:
\begin{equation}
  \label{eq:complex-variables}
  w=b+i r^2 \quad \text{and} \quad  z=x+(b+i r^2)y = x+wy \,.
\end{equation}
The generic solution of two differential operators~\eqref{eq:Cauchy-first}--\eqref{eq:Cauchy-second} is:
\begin{align}
  \label{eq:first-order-generic}
  [\oper{G}f_2](x,y,b,r) & = \sqrt{r}\, \rme^{-\pi \myhbar x^2/(r^2-ib)}\, f_2\left(x+(b+i r^2)y, b+i r^2\right)\\
  \nonumber 
   & = \sqrt{r}\, \rme^{-\pi i \myhbar x^2/w}\, f_2 (z, w)\,,
\end{align}
where \(f_2\) is a holomorphic function of two complex variables. Clearly, with \(b=0\) and \(r=1\) the function \([\oper{G}f_2](x,y,1,0)\) is an element of the (pre-)FSB space on \(\Space{C}{}\).

Furthermore, the integral kernel and metamorphism image space are annihilated by second-order  differential operators:
\begin{align}
  \label{eq:strctural-1}
 \oper{S}_1 &=r^{2}( 4 \pi i   \myhbar \partial_{b} -  \partial_{xx}^{2}) \,;\\
  \label{eq:strctural-2}
 \oper{S}_2 &= -2  \pi i r \myhbar \partial_{r}-  b \partial_{xx}^{2}+\partial_{xy}^{2} - 2  \pi i x \myhbar \partial_{x}-i  \pi  \myhbar I\,.
\end{align}
Of course, the list of annihilators is not exhausted and the above conditions are not independent.
If \( \oper{S}_1 [\oper{G}f_2]=0\) for the generic solution~\eqref{eq:first-order-generic} then the function \(f_2\) has to satisfy the second-order differential equation:
\begin{displaymath}
  w\partial_{zz}^2 f_2(z,w)-4\pi i h z\partial_zf_2(z,w)-4\pi i h w\partial_wf_2(z,w) -2\pi i h f_2(z,w) =0\,.
\end{displaymath}
Equivalently:
\begin{equation}
  \label{eq:structural-generic}
  \partial_wf_2 = \frac{1}{4\pi i h }\partial_{zz}^2 f_2
  - \frac{z}{w} \partial_zf_2 - \frac{1}{2 w} f_2 \,.
\end{equation}
This equation can be reduced to the standard Schr\"odinger equation of a free particle on the line by the change of variables~\amscite{PolyaninNazaikinskii16a}*{\S3.8.3.4}:
\begin{equation}
  \label{eq:structural-change-var}
  (z,w,f_2)\rightarrow \left (\frac{z}{ w}, \frac{1}{w}, \frac{1}{\sqrt{w}} f_2\right) .
\end{equation}
The same equation~\eqref{eq:structural-generic} appears from the condition \( \oper{S}_2 [\oper{G}f_2]=0\). 
Thus, we need only operators \(\oper{C}_1\), \(\oper{C}_2\) and \(\oper{S}_1\) to specify \(f_2\) (and therefore \(\oper{G}f_2\)~\eqref{eq:first-order-generic}).
In the following we call~\eqref{eq:structural-generic} the \emph{structural condition}.

We can check the above characterisation of the image space on the metamorphisms from Example~\ref{ex:initial-metamorphisms}.
\begin{example}
  \begin{enumerate}
  \item For the metamorphism of the wave packet \(\metamorph{P}_{\sigma,\lambda}\) in~\eqref{eq:metamorphism-wave-packet} we find the respective form \(\metamorph{P}_{\sigma,\lambda}  = \oper{G}f_2\) in~\eqref{eq:first-order-generic} as follows:
    \doublechecked
\pyginac
sigma=realsymbol("sigma")
u=realsymbol("u")
h=realsymbol("h")
y=realsymbol("y")
x=realsymbol("x")
r=realsymbol("r")
b=realsymbol("b")
z=symbol("z")
w=symbol("w")
la=realsymbol("lambda")
W=-Pi*sigma*u**2-2*Pi*I*la*u
E=-Pi*sigma*(u-y)**2-2*Pi*h*I*(u-y)*(la-I*sigma*y)/h
# print((W-E).expand().normal())
# - \pi y^{2} \sigma-{(2 i)}  \pi y \lambda
P=W-E-Pi*h*(x+(la-I*sigma*y)/h)**2/(r**2+sigma/h-I*b)+Pi*I*h*x**2/w
PC=P.subs(b==w-I*r**2).expand().normal().subs(x==z-y*w).expand().normal()
print(PC)
\endpyginac
\begin{align}
      \nonumber 
      f_2(z,w)
      &=\frac{1}{\sqrt{- i w+\sigma /\myhbar}} \exp\left(
        \pi \frac{(\lambda ^2 w     +2 z \lambda  w \myhbar      -i z^2 \sigma \myhbar )}{(i \myhbar  w - \sigma ) w}
        \right)\\
      &=\frac{1}{\sqrt{- i w+\sigma /\myhbar}} \exp\left(
       \frac{\pi}{i \myhbar   - \sigma/w}  \left(\lambda ^2 \frac{1}{w}     +2  \lambda   \myhbar  \frac{z}{w}     -i \sigma \myhbar   \frac{z^2}{w^2}\right)
      \right).
      \label{eq:metamorphism-wave-packet-complex}
    \end{align}
  \item The metamorphism of wave~\eqref{eq:wave-transformed} has the representation  \(\metamorph{E}_k=\oper{G}f_2\) for
    \doublechecked
  \begin{equation}
    \label{eq:wave-complex-form}
    f_2(z,w) = \frac{1}{\sqrt{- i w}} \exp\left(-i   k \frac{z }{ w}
      + \frac{  k^{2}}{ 4 \pi i \myhbar}\frac{1}{ w}\right) .
  \end{equation}
  Here, variables~\eqref{eq:structural-change-var} are substituted into the particular solution \(\rme^{\pi i (k^2 t -2k v)}\) of the operator \((4\pi i \myhbar)\partial_t- \partial_{vv}^2\).  Also it is a particular case of~\eqref{eq:metamorphism-wave-packet-complex} for \(\sigma=0\) and \(\lambda=k/(2\pi)\).  
  \end{enumerate}
\end{example}

\subsection{Multidimensional metamorphism}
\label{sec:mult-metam}

It is straightforward to extend metamorphism for functions \(f(u_1,\ldots,u_n)\) on \(\Space{R}{n}\), making \(n\) copies of metamorphism~\eqref{eq:metamorphism} in each variable \(u_1\). \ldots \(u_n\) and receive the function 
on \((\Space[+]{R}{4})^n\). As in the case of the FSB transform, such tensorial approach allow us to immediately extend properties of the one-dimensional metamorphism to an arbitrary finite dimension \(n\).
In particular, for two dimensions the function in the image space shall have the form, cf.~\eqref{eq:first-order-generic}
\begin{equation}
  \label{eq:first-order-generic-doubled}
  [\oper{G}f_4](x_1,y_1,b_1,r_1;x_2,y_2,b_2,r_2) = \sqrt{r_1r_2}\, \rme^{\pi i \myhbar (x_1^2/{w}_1+x_2^2/{w}_2)} f_4\left(z_1,w_1;z_2,w_2\right) ,
\end{equation}
for some function \(f_4\) holomorphic in four complex variables cf.~\eqref{eq:complex-variables}:
\begin{displaymath}
  w_k=b_k+i r_k^2 \quad \text{and} \quad  z_k=x_k+(b_k+i r_k^2)y_k = x_k+w_ky_k \qquad \text{for } k=1, 2.
\end{displaymath}
Furthermore, \(f_4\) needs to satisfy to the structural condition~\eqref{eq:structural-generic} in each pair \((z_1,w_1)\) and \((z_2,w_2)\) of its variables.

Thus, to reduce technicalities we stick in the following to the one-dimensional metamorphism whenever possible.

\begin{example}
  For the plane wave \(E_{k_1 k_2}(u)=\rme^{- i (k_1 u_1+k_2 u_2)}\)  we have, cf.~\eqref{eq:wave-complex-form}:
  \begin{align}
    \nonumber 
    \lefteqn{\metamorph{E}_{k_1k_2} (x_1,y_1,b_1,r_1;x_2,y_2,b_2,r_2)
    = \sqrt{r_1 r_2}\,
      \rme^{-\pi i \myhbar (x_1^2/w_1+x_2^2/w_2)}}\\
    \label{eq:wave-transformed-two}
    & \qquad \qquad \times \frac{1}{\sqrt{- w_1 w_2}}\, 
      \exp\left(-i   k_1 \frac{z_1 }{ w_1} -i   k_2 \frac{z_2 }{ w_2}
      + \frac{  1}{ 4 \pi i \myhbar}\left(\frac{k_1^{2}}{ w_1} + \frac{k_2^{2}}{ w_2}\right)\right) ,
  \end{align}
  where the last line represents the function \(f_4\) in terms of equation~\eqref{eq:first-order-generic-doubled}.
\end{example}

\section{Metamorphism of differential operators}
\label{sec:transm-diff-oper}

We proceed with an application of the metamorphism to differential equations.

\subsection{Reduction of an order of derivations}
\label{sec:reduc-order-second}

A simple calculation shows the intertwining property  for the derivative:
\begin{displaymath}
  (f')\metamorph{\phantom{f}} = \partial_y \metamorph{f}\qquad \text{and thus} \qquad
    (f'')\metamorph{\phantom{f}} = \partial_{yy}^2 \metamorph{f}\,.
\end{displaymath}

Using the annihilation operators~\eqref{eq:Cauchy-first} and~\eqref{eq:strctural-1} we can reduce the second derivative to the first order operator using the expression
\begin{align*}
  \partial_{yy}^2  = \oper{D}_0 - ir \left((b+ir^{2})\partial_x + \partial_y\right)\oper{C}_1 - \frac{1}{r^2}(b+ir^{2})^{2}\oper{S}_1\,,
\end{align*}
where the first order differential operator \(\oper{D}_0\) is
\begin{displaymath}
\oper{D}_0 =  2 \pi i\myhbar \left(
   (b+   i r^{2}) x \gpartial_{0} 
+    x   \gpartial_{1}
  +2 (b+   i r^{2})^2 \gpartial_{2}
    +{(b+   i r^{2})} I
\right) .
\end{displaymath}
Thus we obtain the following lemma
\begin{lem}
  Metamorphism~\eqref{eq:metamorphism} intertwines the second-order derivative to the first-order differential operator:
  \begin{displaymath}
    (f'')\metamorph{\phantom{f}} = \oper{D}_0 \metamorph{f}.
  \end{displaymath}
\end{lem}
The operator \(\oper{D}_0\) becomes explicitly transparent if we use the general solution \([\oper{G}f_2]\) in~\eqref{eq:structural-generic}.
\begin{prop}
  \label{pr:second-transmutes-general}
  For a function \(f(u)\), let its metamorphism~\eqref{eq:metamorphism} be \(\metamorph{f}=\oper{G} f_2\) of~\eqref{eq:structural-generic} for some holomorphic function \(f_2\) of two complex variables. Then:
  \begin{displaymath}
    (f'')\metamorph{\phantom{f}}=\oper{G} (\oper{D} f_2)\,,
  \end{displaymath}
  where for \(f_2(z,w)\):
  \begin{equation}
    \label{eq:second-metamorphism-generic}
    \oper{D} f_2 = 2 \pi i  \myhbar w {(2  z \partial_{z}f_2 +2  w \partial_{w}f_2+f_2)}\,.
  \end{equation}
\end{prop}

\subsection{Application to the Helmholtz equation}
\label{sec:appl-helmh-equat}

As mentioned in  \S~\ref{sec:mult-metam}, for a function \(f(u_1,u_2)\) on \(\Space{R}{2}\) we can repeatedly apply metamorphism~\eqref{eq:metamorphism} in each variable \(u_1\) and \(u_2\). All the above calculations remain valid for the doubled set of variables, of course.
The straightforward application of Prop.~\ref{pr:second-transmutes-general} implies the following.
\begin{cor}
  \label{co:helmholtz-2D}
  Let \(f(u_1,u_2)\) be a solution of the Helmholtz equation:
  \begin{equation}
    \label{eq:Helmholtz2}
    \partial_1^2 f+ \partial_2^2 f - k^2 f = 0.
  \end{equation}
  Then \(\metamorph{f}=\oper{G}f_4\) for a function \(f_4\) holomorphic in four complex variables, which satisfies the first-order differential equation:
  \begin{equation}
    \label{eq:helmholtz-transmuted}
    2w_1 {(  z_1 \partial_{z_1}f_4 +  w_1 \partial_{w_1}f_4)}
    +  2w_2 {(  z_2 \partial_{z_2}f_4 +  w_2 \partial_{w_2}f_4)}
    +\left(w_1 + w_2 + \frac{k^2}{2 \pi i  \myhbar}\right) f_4=0\,.
  \end{equation}
\end{cor}

The presence of the Euler operators for pairs  \((z_1,w_1)\) and \((z_2,w_2)\) in~\eqref{eq:helmholtz-transmuted} suggests to look for a solution of~\eqref{eq:helmholtz-transmuted} in the form:
\begin{displaymath}
  f_4(z_1,w_1; z_2 w_2) = \phi\left(w_1,w_2, \frac{z_1}{w_1},\frac{z_2}{w_2}\right).
\end{displaymath}
Then,  equation~\eqref{eq:helmholtz-transmuted} reduces to the following differential equation for \(\phi\):
\begin{displaymath}
  2w_1^2\, \partial_1 \phi + 2w_2^2\, \partial_2 \phi + \left(w_1 + w_2 + \frac{k^2}{2 \pi i  \myhbar} \right) \phi =0 \,.
\end{displaymath}
The last equation gives the following generic solution of~\eqref{eq:helmholtz-transmuted}, cf.~\amscite{PolyaninZaitsevMoussiaux02}*{\S4.8.2.4}:
\begin{equation}
  \label{eq:solution-helmholtz-reduced}
  f_4(z_1,w_1;z_2,w_2) = \frac{\rme^{k^2(w_1^{-1}+w_2^{-1})/(8 \pi i  \myhbar )}}{\sqrt{w_1 w_2}}  \, f_3\left(\frac{1}{w_1}-\frac{1}{w_2}, \frac{z_1}{w_1},\frac{z_2}{w_2}\right) ,
\end{equation}
with a holomorphic function \(f_3\) of three complex variables.   Taking in account~\eqref{eq:first-order-generic-doubled} the full metamorphism is
\begin{align}
  \nonumber 
  \metamorph{f}(z_1,w_1;z_2,w_2)
  &=
      \sqrt{\frac{r_1 r_2}{{w_1 w_2}}}\, 
    \exp\left({-\pi i \myhbar \left(\left(x_1^2+ \frac{\kbar^2}{2\myhbar^2}\right)\frac{1}{w_1}
    +\left(x_2^2+\frac{\kbar^2}{2\myhbar^2}\right)\frac{1}{w_2}\right)}\right)\\
  \label{eq:solution-helmholtz-full}
  & \qquad \times 
    f_3\left(\frac{1}{w_1}-\frac{1}{w_2}, \frac{z_1}{w_1},\frac{z_2}{w_2}\right).
\end{align}

Finally, we check the structural conditions~\eqref{eq:strctural-1} on the last expression. An application of the operator \(\oper{S}_1\) in variables \((z_1,w_1)\) to \(\oper{G}f_4\) from~\eqref{eq:solution-helmholtz-reduced} produces the equation:
\begin{equation}
  \label{eq:structural-1-helmholtz}
  4 \pi i \myhbar   \partial_{1} f_3 -   \partial_{22}^{2}f_3 - \half  k^{2}   f_3 = 0 \,.
\end{equation}
Similarly an application of the operator \(\oper{S}_1\) in variables \((z_2,w_2)\) produces:
\begin{equation}
  \label{eq:structural-2-helmholtz}
4 \pi i \myhbar   \partial_{1} f_3 + \partial_{33}^{2}f_3 + \half  k^{2}   f_3 = 0 \,.
\end{equation}
These are Schr\"odinger equations of a free particle with one degree of freedom. Note the opposite flow of time in them.

Also, a direct calculation shows that an application of the Helmholtz operator to \(\oper{G} f_4\) for some \(f_4\) from~\eqref{eq:solution-helmholtz-reduced} reduces to
\begin{displaymath}
  \partial_{22}^{2}f_3 + \partial_{33}^{2}f_3 +   k^{2}   f_3 = 0\,,
\end{displaymath}
which is the difference of~\eqref{eq:structural-1-helmholtz} and \eqref{eq:structural-2-helmholtz}, and thus follows from them.
\begin{thm}
  Let \(f(u_1,u_2)\) be a solution to the Helmholtz equation~\eqref{eq:Helmholtz2}, then \(\metamorph{f}=\oper{G} f_4\), where \(f_4\) is a function of the form~\eqref{eq:solution-helmholtz-reduced} for some holomorphic \(f_3\), which  satisfies two structural conditions~\eqref{eq:structural-1-helmholtz}--\eqref{eq:structural-2-helmholtz}.
\end{thm}

\begin{figure}
  \centering
    \includegraphics[scale=1]{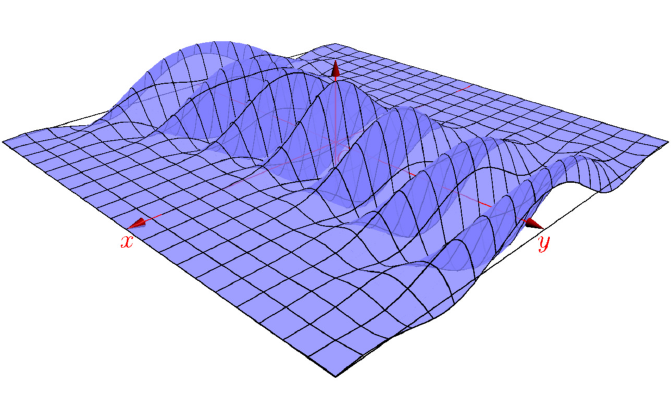}
    \includegraphics[scale=1]{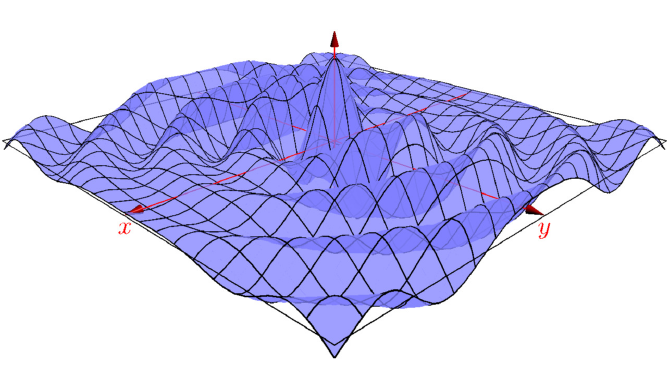}
    \caption{Two Gaussian beams: narrow and wide.}
  \label{fig:two-gaussian-beams}
\end{figure}

Therefore, one can construct a particular solution of the Helmholtz equation in the following steps:
\begin{enumerate}
\item Take some particular solutions of the Schr\"odinger equation and blend them into a joint solution \(f_3\) to the system of~\eqref{eq:structural-1-helmholtz}--\eqref{eq:structural-2-helmholtz}.
\item Use the obtained function \(f_3\) to construct solution \(f_4\) from~\eqref{eq:solution-helmholtz-reduced} of the transmuted Helmholtz equation~\eqref{eq:helmholtz-transmuted}.
\item Use~\eqref{eq:first-order-generic-doubled} to re-create the full solution \(\oper{G}f_4\) on \(\Space{R}{8}\) from the above \(f_4\),
\item Then, either
  \begin{itemize}
  \item employ \(\oper{G}f_4\) to analyse the problem in terms of
    coordinates on \(\Space{R}{8}\); or
  \item use integral transform~\eqref{eq:metamorphism-inverse} to recover the solution in terms of original coordinates \((u_1,u_2)\).
\end{itemize}
\end{enumerate}
\begin{example}
  Many partial solutions of the heat/Schr\"odinger equation~\amscite{PolyaninNazaikinskii16a}*{\S~3.1.3} are variations of two main themes: the plane wave-type and the fundamental solution (the Gaussian wave packet). We will treat both of them now.
  \begin{enumerate}
  \item  Let us start from the partial solutions~\amscite{PolyaninNazaikinskii16a}*{\S~3.1.3}
    \begin{equation}
      \label{eq:plane-wave-partial}
      \phi_j(z_j, w_j)=\rme^{(k_j^2-k^2/2)/(4\pi i \myhbar)/w_j-  i k_j z_j/w_j}\,,\qquad\text{for } j=1,2\,.
    \end{equation}
    of a plane wave-type of the structural equations~\eqref{eq:structural-1-helmholtz}--~\eqref{eq:structural-2-helmholtz}.  Then we have a representation 
    \begin{displaymath}
      \phi_1(z_1, w_1) \phi_2(z_2, w_2) = f_3\left(\frac{1}{w_2}-\frac{1}{w_1}, \frac{z_1}{w_1},\frac{z_2}{w_2}\right)
    \end{displaymath}
    if and only if \(\partial_2\phi_1= -\partial_2\phi_2\), that is \(k_1^2+k_2^2=k^2\). Thus,
    \begin{align}
      \nonumber 
      f_4(z_1,w_1;z_2,w_2)&= \frac{\rme^{(k_1^2+k_2^2)(w_1^{-1}+w_2^{-1})/(8 \pi i  \myhbar )}}{\sqrt{w_1 w_2}}\\
      \nonumber 
      &\qquad\times
      \exp\left(\frac{k_1^2-k_2^2}{8\pi i \myhbar}\left(\frac{1}{w_1}-\frac{1}{w_2}\right)
        - i k_1 \frac{z_1}{w_1} - i k_2 \frac{z_2}{w_2}\right)\\
      \label{eq:plane-wave-complex-deduced}      
      &= \frac{1}{\sqrt{w_1 w_2}}\,
        \exp\left(\frac{k_1^2}{4\pi i \myhbar}\frac{1}{w_1}-i k_1 \frac{z_1}{w_1}\right)
        \exp\left(\frac{k_2^2}{4\pi i \myhbar}\frac{1}{w_2}-i k_2 \frac{z_2}{w_2}\right) .
    \end{align}
    In other words, we have recovered the metamorphism of a plane wave in the plane \((u_1,u_2)\) along the direction \((k_1,k_2)\)~\eqref{eq:wave-transformed-two}. 
  \item Consider the fundamental solutions~\amscite{PolyaninNazaikinskii16a}*{\S~3.1.3}
    \begin{displaymath}
      \psi_j(z_j, w_j)= \sqrt{\myhbar w_j} \exp\left(\pi i \myhbar \frac{z_j^2}{w_j}+\frac{k^2}{8\pi i \myhbar}\frac{1}{w_j}\right),\qquad\text{for } j=1,2\,.
    \end{displaymath}
    of the structural conditions~\eqref{eq:structural-1-helmholtz}--~\eqref{eq:structural-2-helmholtz}. We note that each \(\psi_j\) is a superposition of the plane waves~\eqref{eq:plane-wave-partial} over the wavenumber  \(k_j\) with a Gaussian density. Thus, we can build the respective function as similar superposition of the plane waves~\eqref{eq:plane-wave-complex-deduced} considering the allowed range \(k_1\in [-k, k]\) as follows:
    \begin{align*}
      \lefteqn{f_3\left(\frac{1}{w_2}-\frac{1}{w_1}, \frac{z_1}{w_1},\frac{z_2}{w_2}\right)}\\
      &\qquad=  \int\limits_{-k}^k  \exp\left(\frac{2k_1^2-k^2}{8\pi i \myhbar}\left(\frac{1}{w_1}-\frac{1}{w_2}\right) - i k_1 \frac{z_1}{w_1} \mp i \sqrt{k^2-k_1^2} \frac{z_2}{w_2}\right) \rme^{-a k_1^2}  \,\rmd k_1\,.
    \end{align*}
    Illustration of this beam for two different values of the parameter \(a\) is given in Fig.~\ref{fig:two-gaussian-beams}.
  \end{enumerate}
\end{example}

\subsection{Helmholtz operator in three dimensions}
\label{sec:helmh-oper-three}

Solution of the Helmholtz equation in higher dimensions follows  the same steps which were used in two dimensions.
Using the three copies of transformation~\eqref{eq:metamorphism} we obtain the following analog of Cor.~\ref{co:helmholtz-2D}.

\begin{cor}
  Let \(f(u_1,u_2)\) be a solution of the Helmholtz equation:
  \begin{displaymath}
    \partial_1^2 f + \partial_2^2 f + \partial_3^2 f - k^2 f = 0.
  \end{displaymath}
  Then \(\metamorph{f}=\oper{G}f_6\) for a function \(f_6\) holomorphic in four variables which satisfies the first-order differential equation:
  \begin{align}
    \nonumber 
    2w_1 {(  z_1 \partial_{z_1}f_6 +  w_1 \partial_{w_1}f_6)}
    +  2w_2 {(  z_2 \partial_{z_2}f_6 +  w_2 \partial_{w_2}f_6)}
    +  2w_3 {(  z_3 \partial_{z_3}f_6 +  w_3 \partial_{w_3}f_6)}&\\
    \label{eq:helmholtz3-transmuted}
    +\left(w_1 + w_2 + w_3 + \frac{k^2}{2 \pi i  \myhbar}\right) f_6&=0\,.
  \end{align}
\end{cor}

As before, the presence of the Euler operators for pairs  \((z_j,w_j)\), \(j=1, 2,3\) in~\eqref{eq:helmholtz3-transmuted} suggests to look for a solution of~\eqref{eq:helmholtz3-transmuted} of the form:
\begin{displaymath}
  f_6(z_1,w_1; z_2 w_2; z_3 w_3) = \phi\left(w_1,w_2, w_3, \frac{z_1}{w_1}, \frac{z_2}{w_2}, \frac{z_3}{w_3}\right).
\end{displaymath}
Then,  equation~\eqref{eq:helmholtz3-transmuted} reduces to the following differential equation for \(\phi\):
\begin{displaymath}
  2w_1^2\, \partial_1 \phi + 2w_2^2\, \partial_2 \phi + 2w_3^2\, \partial_3 \phi + \left(w_1 + w_2 + w_3 + \frac{k^2}{2 \pi i  \myhbar} \right) \phi =0 \,.
\end{displaymath}
It has the generic solution~\amscite{PolyaninZaitsevMoussiaux02}*{\S~8.8.2.1} of the form:
\begin{align}
  \nonumber 
  f_6(z_1,w_1;z_2,w_2;z_3,w_3) &= \frac{\rme^{k^2(w_1^{-1}+w_2^{-1}+w_3^{-1})/(12 \pi i  \myhbar )}}{\sqrt{w_1 w_2 w_3}}  \\
  \label{eq:solution-helmholtz3-reduced}
  & \qquad \times f_5\left(\frac{1}{w_1}-\frac{1}{w_3}, \frac{1}{w_2}-\frac{1}{w_3}, \frac{z_1}{w_1}, \frac{z_2}{w_2}, \frac{z_3}{w_3}\right) ,
\end{align}
with a holomorphic function \(f_5\) of five complex variables. Furthermore, we again need to satisfy the following three structural conditions:
\begin{align}
    \label{eq:structural-1-helmholtz-3D}
4 \pi i \myhbar   \partial_{1} f_5 +   \partial_{33}^{2}f_5 + \third  k^{2}   f_5 &= 0 \,,\\
    \label{eq:structural-2-helmholtz-3D}
4 \pi i \myhbar   \partial_{2} f_5 +   \partial_{44}^{2}f_5 + \third  k^{2}   f_5 &= 0 \,,\\
    \label{eq:structural-3-helmholtz-3D}
4 \pi i \myhbar   \partial_{1} f_5 + 4 \pi i \myhbar   \partial_{2} f_5 -   \partial_{55}^{2}f_5 - \third  k^{2}   f_5 &= 0 \,.
\end{align}
A similar form of solutions for the Helmholtz equation can be obtained in any dimension.

\section{Discussion and conclusions}
\label{sec:conclusions}

We have presented integral transform~\eqref{eq:metamorphism} acting \(\FSpace{L}{2}(\Space{R}{}) \rightarrow \FSpace{L}{2}(\Space[+]{R}{4})\). The image space  of the metamorphism consists of functions satisfying  certain holomorphic conditions~\eqref{eq:strctural-1}--\eqref{eq:strctural-2} and Schr\"odinger-type structural equations~\eqref{eq:strctural-1} (which is equivalent to~\eqref{eq:strctural-2}). These restrictions imply that the metamorphism intertwines the second derivative with the first order differential operator~\eqref{eq:second-metamorphism-generic}. Reduction of the order by \(1\) replaces a second-order PDE by a first-order equation (in higher dimensions), which can be solved by the method of characteristics. These techniques applied to the Helmholtz equations produce the metamorphism of a generic solution in the form~\eqref{eq:solution-helmholtz-full}.

The metamorphism is a realisation of the coherent state transform~\amscite{AliAntGaz14a}*{Ch.~8} for a unitary irreducible representation of the semi-direct product of the Heisenberg group \(\Space{H}{1}\) and the \(ax+b\) group acting  on  \(\Space{H}{1}\) by symplectomorphisms. Restricting the representation to various subgroups we obtain many familiar integral transforms, cf. Rem.~\ref{re:generalisations}. Their combined power is naturally inherited by the metamorphism. The group-theoretic aspects of metamorphism will be presented elsewhere. 

\emph{Summing up:} metamorphism presents a new general method of solving and studying various PDEs and deserves further thoughtful investigation.

\section*{Acknowledgments}
\label{sec:acknowledgments}
I am grateful to Dr.~A.V.~Kisil for fruitful collaboration on this topic.
Prof.~S.M.~Sitnik provided enlightening information on the Gauss--Fresnel integral.
An anonymous referee provided many useful comments.

\small
\providecommand{\noopsort}[1]{} \providecommand{\printfirst}[2]{#1}
  \providecommand{\singleletter}[1]{#1} \providecommand{\switchargs}[2]{#2#1}
  \providecommand{\irm}{\textup{I}} \providecommand{\iirm}{\textup{II}}
  \providecommand{\vrm}{\textup{V}} \providecommand{\cprime}{'}
  \providecommand{\eprint}[2]{\texttt{#2}}
  \providecommand{\myeprint}[2]{\texttt{#2}}
  \providecommand{\arXiv}[1]{\myeprint{http://arXiv.org/abs/#1}{arXiv:#1}}
  \providecommand{\doi}[1]{\href{http://dx.doi.org/#1}{doi:
  #1}}\providecommand{\CPP}{\texttt{C++}}
  \providecommand{\NoWEB}{\texttt{noweb}}
  \providecommand{\MetaPost}{\texttt{Meta}\-\texttt{Post}}
  \providecommand{\GiNaC}{\textsf{GiNaC}}
  \providecommand{\pyGiNaC}{\textsf{pyGiNaC}}
  \providecommand{\Asymptote}{\texttt{Asymptote}}

\end{document}